\documentclass[11pt, a4paper]{article}
\usepackage{fullpage}
\usepackage{amsfonts}
\usepackage{amssymb}
\usepackage{amsmath}
\usepackage{amsthm}
\usepackage{graphicx}
\usepackage{bm}
\usepackage[makeroom]{cancel}
\usepackage{enumitem}
\usepackage{url}
\usepackage[margin=.9in]{geometry}
\usepackage{amsopn}
\usepackage{mathtools}
\usepackage{hyperref}
\usepackage{doi}
\usepackage{cite}
\usepackage{overpic}
\usepackage[symbol]{footmisc}

\usepackage{listings}
\usepackage{xcolor}
\definecolor{codegreen}{rgb}{0,0.6,0}
\definecolor{codegray}{rgb}{0.5,0.5,0.5}
\definecolor{codepurple}{rgb}{0.58,0,0.82}
\definecolor{backcolour}{rgb}{0.95,0.95,0.92}
\lstdefinestyle{mystyle}{
    backgroundcolor=\color{backcolour},   
    commentstyle=\color{codegreen},
    keywordstyle=\color{magenta},
    numberstyle=\tiny\color{codegray},
    stringstyle=\color{codepurple},
    basicstyle=\ttfamily\normalsize,
    breakatwhitespace=false,         
    breaklines=true,                 
    captionpos=b,                    
    keepspaces=true,                 
    numbers=left,                    
    numbersep=3pt,                  
    showspaces=false,                
    showstringspaces=false,
    showtabs=false,                  
    tabsize=4,
}
\usepackage{amsmath}

\DeclareMathOperator{\sech}{sech}
\DeclareMathOperator{\arctantwo}{atan2}

\begin{document}

\begin{center}
    \Large \bf Phasor notation of Dynamic Mode Decomposition
\end{center}
\begin{center}
    Karl Lapo$^{1}$\footnote[1]{Corresponding authors (karl-eric.lapo@uibk.ac.at)},
    Samuele Mosso$^{1}$,
    J. Nathan Kutz$^{2,3}$
\end{center}
\begin{center}
    \scriptsize{
    ${}^1$ Department of Cyrospheric and Atmospheric Sciences, University of Innsbruck, Innsbruck, Austria \\ 
    ${}^2$ Department of Electrical and Computer Engineering, University of Washington, Seattle, WA 98195, United States \\
     ${}^3$ Department of Applied Mathematics, University of Washington, Seattle, WA 98195, United States     }
\end{center}


\begin{abstract}
Dynamic Mode Decomposition (DMD) is a powerful, data-driven method for diagnosing complex dynamics. Various DMD algorithms allow one to fit data with a low-rank model that decomposes it into a sum of coherent spatiotemporal patterns. Nominally, each rank of the DMD model is interpreted as a complex, stationary spatial mode modulated by a single set of complex time dynamics (consisting of exponential growth/decay and oscillation), and an amplitude. However, the specifics of how these DMD components are interpreted do not appear to be consistent with the information actually present in the DMD decomposition or the underlying data. While there is a clear physical interpretation for the complex time dynamics, there is practically no guidance on the complex spatial modes. To resolve these issues, we introduce the phasor notation of the DMD model for conjugate pair DMD modes, which results in a strictly positive and real spatial pattern as well as spatiotemporal waveform. The phasor notation terms result in an interpretable DMD model that provides a more complete diagnoses of the model components, as demonstrated on a toy model. This DMD interpretation needs to be adjusted for DMD variants which alter the relationship between the DMD model and the data, such as those that window data in time. We derive the phasor notation terms for one such method, multi-resolution Coherent Spatiotemporal Scale-separation, and demonstrate the new terms by interpreting a multi-scale data set.
\end{abstract}


\section{Introduction}
There can be a substantial ``interpretability gap" between the guidance provided for the use of a mathematical method and the reality of its use on real-world data.  Guidance on the use of a method is usually provided in the form of applications using simple or toy data. This guidance then informs the application of the method in more complicated settings such as in the physical sciences. However, the realities of physical data require a physical interpretation of the new method that is somewhat distinct from the initial guidance. A famous example of closing the interpretability gap is the introduction of continuous wavelets transforms (CWT) to the earth sciences community by Torrence and Compo\cite{Torrence1998}, in which methods for quantitative, as opposed to qualitative, analysis with clear physical interpretations were provided.  With the rapid growth of complex data prevalent in modern science and engineering, there is a substantial need for methods to be more clearly understood and interpreted in the context of such complexity. 

In this work we close this interpretability gap for Dynamic Mode Decomposition (DMD). DMD is a data-driven model that seeks to find spatial modes with shared time dynamics prescribed by a complex exponential. Historically, the definition of DMD hinges on the view of DMD as an algorithm, as first espoused in \cite{Schmid2010_dmd}. However, several solvers now exist for the DMD model which do not rely on the original algorithm e.g., \cite{Askham2018, ichinagaPyDMDPythonPackage2024a}. As such, we argue that DMD should not be viewed as an algorithm, but as a model for which multiple algorithms can find reasonable solutions. Clearly distinguishing between the DMD solver and model allows one to more directly elucidate the physical interpretation of DMD.

The DMD model fits data with a linear superposition of the following form
\begin{equation}\label{eq:DMD-model}
\tilde{\mathbf{x}}(t)=\sum^{r}_{j=1}\boldsymbol{\phi}_j e^{\lambda_jt}b_j
\end{equation}
where $\tilde{\mathbf{x}}$ is a low-rank approximation of the original snapshot matrix with components $\mathbf{x}$ which are collected across time $t$. The summation is performed over complex eigenvalues ($\lambda_j$) and eigenvectors ($\boldsymbol{\phi}_j$) which are scaled by an amplitude, $b_j$. The DMD model provides an approximation of rank $r$ of the original data while also providing an interpretable model of the linearized dynamics. Typically, it is stated that $\boldsymbol{\phi}_j$, the eigenvector, is a coherent spatial mode governed by a single set of time dynamics given by $\lambda_j$, the eigenvalue, with the amplitude, $b_j$. Note that we denote vectors in space using bold symbols.

For most DMD fits on real-valued data, the DMD modes often consist of conjugate pairs such that there are only $\frac{r}{2}$ unique eigenvector/eigenvalue pairs, or modes. There can often be a DC component ($\lambda_j=0$) as well that is unpaired in the DMD reduction. For the time dynamics, the imaginary component of the eigenvalue, $Im(\lambda)$, yields the oscillatory component of the time dynamics while the real part, $Re(\lambda)$, gives the exponential growth/decay. For the spatial modes, $\boldsymbol{\phi}_j$ is also complex. A typical interpretation is that the spatial mode modulated by the time dynamics is just the real component, $Re(\phi)$, while $Im(\phi)$ is neglected e.g., as implemented in \cite{ichinagaPyDMDPythonPackage2024a} or \cite{Kutz_DMD-textoobk_2016}. Finally, $b$ should be a real, positive value.  However, this parsimonious interpretation fails to correctly interpret the complex components of $\boldsymbol{\phi}$, nor what these components mean in combination, especially in the context of the conjugate pairs of $\boldsymbol{\phi}$ \cite{Kutz_DMD-textoobk_2016, Taira2017, Taira2020}. These problems become especially apparent when trying to interpret a collection of DMD models which vary in time, such as from multi-resolution DMD (mrDMD) \cite{Dylewsky2019_dmd}, non-stationary DMD (nsDMD) \cite{ferreNonStationaryDynamicMode2023c}, and multi-resolution Coherent Spatiotemporal Scale-separation (mrCOSTS) \cite{lapoMethodUnsupervisedLearning2025}. In these frameworks, we have discovered various behaviors which are challenging to reconcile with this commonly asserted DMD interpretation, such as $Re(\phi)$ changing sign between adjacent time windows despite constant physics.

In Sect. \ref{sect:dmd-phasor-notation} we introduce and formalize a physically-consistent way of understanding the DMD model using a phasor notation of the DMD model. Using DMD in phasor form we recast our intuitive understanding of the interpretable DMD model, enabling the direct recovery of pieces of information that were asserted to exist in the DMD model but were in practice difficult to recover, as well as other neglected pieces of information present in the DMD model that were masked by the `classical' DMD interpretation. The phasor notation interpretation is demonstrated in Sect. \ref{sect:uniscale-results} using a simple, uniscale toy model. In the second part of the study, Sect. \ref{sect:mrCOSTS-phasor-notation}, we show how the the phasor notation of the DMD model needs to be adapted for DMD variants which change the definition of the DMD model, e.g. by changing the DMD model terms with time. For this demonstration we use mrCOSTS) \cite{lapoMethodUnsupervisedLearning2025}, a windowed, hierarchical decomposition based on DMD, in order to demonstrate how to adapt this interpretation for these DMD variants. The derived expressions are then demonstrated on a multi-scale toy example.

\section{The DMD model in phasor notation}\label{sect:dmd-phasor-notation}

To derive the phasor representation of DMD, we need to represent the eigenvalue/eigenvector pairs as complex objects which we notate as 
\begin{equation}
    \boldsymbol{\phi}=\boldsymbol{\phi}^R+i\boldsymbol{\phi}^I
\end{equation} 
and $\lambda=\mu+i\omega$ for simplicity. The DMD model Eq. \ref{eq:DMD-model} can then be re-written as
\begin{equation}\label{eq:DMD-complex-exponential}
\tilde{\mathbf{x}}(t)=\sum^{r}_{j=1} b_j \left( \boldsymbol{\phi}^R_j + i\boldsymbol{\phi}^I_j \right) e^{(\mu_j + i\omega_j) t}.
\end{equation}
Next, we can use Euler's formula to re-express the imaginary part of the exponent in Eq. \ref{eq:DMD-complex-exponential} as a complex sum of cosine and sine, 
\begin{equation}
\tilde{\mathbf{x}}(t)=\sum^{r}_{j=1} e^{\mu_jt} b_j \left( \boldsymbol{\phi}^R_j + i\boldsymbol{\phi}^I_j \right) \left(\cos{\omega_j t} + i\sin{\omega_j t} \right)
\end{equation}
before expanding the individual terms
\begin{equation}\label{eq:dmd_sum-sine-cosine}
\tilde{\mathbf{x}}(t)=\sum^{r}_{j=1} e^{\mu_jt} b_j \left( \boldsymbol{\phi}^R_j \cos{\omega_j t} + i\boldsymbol{\phi}^I_j\cos{\omega_j t} + i\boldsymbol{\phi}^R_j \sin{\omega_j t} - \boldsymbol{\phi}^I_j \sin{\omega_j t} \right).
\end{equation}
Here it is important to point out that the imaginary components should not be immediately neglected.

\subsection{Considering conjugate pairs}\label{sect:conjugate-pairs} One of the complexities of interpreting the DMD model is the presence of the conjugate pairs. The DMD conjugate pairs refer to both complex components of the DMD model, $\boldsymbol{\phi}$ and $\lambda$. Let us consider a pair of DMD ranks that are conjugate pairs, $j=1$ and $j=2$, denoted as $\boldsymbol{\phi}_2=\boldsymbol{\phi}_1^*=\boldsymbol{\phi}^R_1-i\boldsymbol{\phi}^I_1$ and $\lambda_2=\lambda_1^*=\mu_1-iq_1$. First, let us consider the reconstruction of the data provided by a single DMD rank, $j=1$,
\begin{equation}\label{eq:dmd_pair-1}
\tilde{\mathbf{x}}_{1}(t)= e^{\mu_1t} b_1  \left( \boldsymbol{\phi}^R_1 \cos{\omega_1 t} + i\boldsymbol{\phi}^I_1\cos{\omega_1 t} + i\boldsymbol{\phi}^R_1 \sin{\omega_1 t} - \boldsymbol{\phi}^I_1 \sin{\omega_1 t} \right).
\end{equation}
For the conjugate pair, $j=2$, Eq. \ref{eq:DMD-complex-exponential} becomes
\begin{equation}\label{eq:dmd_pair-2}
\tilde{\mathbf{x}}_{2}(t)= e^{\mu_1t} b_1  \left( \boldsymbol{\phi}^R_1 \cos{\omega_1 t} - i\boldsymbol{\phi}^I_1\cos{\omega_1 t} - i\boldsymbol{\phi}^R_1 \sin{\omega_1 t} - \boldsymbol{\phi}^I_1 \sin{\omega_1 t} \right).
\end{equation}
Adding together the contribution from each complex pair, Eq. \ref{eq:dmd_pair-1} and Eq. \ref{eq:dmd_pair-2}, simply reduces to
\begin{equation}\label{eq:added-pairs}
\tilde{\mathbf{x}}_{1}(t) + \tilde{\mathbf{x}}_{2}(t)= 2 e^{\mu_1t} b_1 \left[\boldsymbol{\phi}^R_1 \cos(\omega_1t) - \boldsymbol{\phi}^I_1 \sin(\omega_1t)\right].
\end{equation}
The imaginary components of the solution cancel out completely and can be safely neglected when considering DMD fits with conjugate pairs. Note that this statement is not equivalent to saying the imaginary components of the DMD model cancel out, since $\boldsymbol{\phi}^I$ contributes to $\tilde{\bf x}$, counter to the intuition and interpretations of standard DMD. Further, $\tilde{\bf x}_1$ will have an imaginary component when considered in isolation from its conjugate pair, $\tilde{\bf x}_2$. Thus, the DMD is only strictly real for solutions that are strictly conjugate pairs. Conversely, in the case that non-conjugate pairs are found, $\tilde{\mathbf{x}}$ can include an imaginary component. The eigenvalue constraints implemented in the PyDMD package \cite{ichinagaPyDMDPythonPackage2024a} allows one to force conjugate pair solutions for the variable projection solver \cite{Askham2018, ichinagaPyDMDPythonPackage2024a}, enabling this requirement to be met robustly. A similar constraint is not implemented for exact DMD solvers to our knowledge. 

\subsection{Phasor notation}\label{sect:phasor-notation}
From this line of reasoning we can express the DMD model more generally as
\begin{equation}\label{eq:DMD-phasor-notation-as-sum}
\tilde{\mathbf{x}}(t)=\sum^{r}_{j=1} b_j e^{\mu_jt} \left( \boldsymbol{\phi}^R_j \cos{\omega_j t} - \boldsymbol{\phi}^I_j \sin{\omega_j t} \right).
\end{equation}
When considering conjugate pair solutions this notation has a result identical to the original DMD notation since the imaginary components of the solution will cancel.

Eq. \ref{eq:DMD-phasor-notation-as-sum} can be further simplified into a single cosine with a phase shift using the trigonometric identify for harmonic addition of sine and cosine
\begin{equation}\label{eq:trig-id_linear-sum}
a \cos(\theta) + b \sin(\theta) = c \cos(\theta - \varphi)
\end{equation}
which has an amplitude of $c = \sqrt{a^2 + b^2}$ and phase shift of $\varphi = \arctantwo(b, a)$ (n.b. that $b$ here refers to the amplitude of the sine wave and that the $\arctantwo$ function handles any sign issues). Using this identity we can re-write Eq. \ref{eq:DMD-phasor-notation-as-sum} as
\begin{equation}
\tilde{\mathbf{x}}(t)=\sum^{r}_{j=1} b_j e^{\mu_jt} \sqrt{\boldsymbol{\phi}^{R,2}_j + \boldsymbol{\phi}^{I,2}_j} \cos{\left(\omega_j t - \arctantwo{(\boldsymbol{\phi}^I_j, \boldsymbol{\phi}^R_j)}\right)}.
\end{equation}
However, in practice with DMD solutions, $\omega$ and $\varphi$ have their signs already flipped relative to each other such that numerical implementations need to be expressed as 
\begin{equation}\label{eq:DMD-phasor-notation}
\tilde{\mathbf{x}}(t)=\sum^{r}_{j=1} b_j e^{\mu_jt} \sqrt{\boldsymbol{\phi}^{R,2}_j + \boldsymbol{\phi}^{I,2}_j} \cos{\left(\omega_j t + \arctantwo{(\boldsymbol{\phi}^I_j, \boldsymbol{\phi}^R_j)}\right)}.
\end{equation}
Finally, the reconstruction of a conjugate pair is notated as
\begin{equation}\label{eq:DMD-phasor-notation_pairs}
\tilde{\mathbf{x}}_{j} + \tilde{\mathbf{x}}_{j+1}(t)=2b_j e^{\mu_jt} \sqrt{\boldsymbol{\phi}^{R,2}_j + \boldsymbol{\phi}^{I,2}_j} \cos{\left(\omega_j t + \arctantwo{(\boldsymbol{\phi}^I_j, \boldsymbol{\phi}^R_j)}\right)}.
\end{equation}

\subsection{The new DMD interpretation} 
This new formulation of the DMD model allows us to extract physically interpretable information with a simple intuitive meaning. The term, $\boldsymbol{\varphi}$, is the spatially-varying phase shift given by
\begin{equation}\label{eq:phase-shift-single-mode}
\boldsymbol{\varphi}_j = \arctantwo(\boldsymbol{\phi}^I_j, \boldsymbol{\phi}^R_j).
\end{equation}
By itself, this term is not novel, for instance its interpretation is specified in \cite{Taira2020} and \cite{proctorDiscoveringDynamicPatterns2015}. However, $\boldsymbol{\varphi}$ is not amenable to a simple, intuitive interpretation on its own due to the angle sum between the components $\boldsymbol{\varphi}$ and $\omega t$. Instead, it must be interpreted as a component of a spatiotemporal waveform, $\boldsymbol{\mathcal{W}}$, given by
\begin{equation}\label{eq:DMD-waveform}
\boldsymbol{\mathcal{W}}_j = \cos(\omega_jt + \boldsymbol{\varphi}_j).   
\end{equation}
which is modulated in space by $\boldsymbol{\varphi}$ and in time by $\omega$. This interpretation is still similar to the original DMD interpretation but with the important clarification that while the waveform of a given mode has a single oscillatory frequency in time it can have a complicated spatiotemporal structure due to the angle addition. We demonstrate in Sect. \ref{sect:uniscale-results} how interpreting $\boldsymbol{\varphi}$ alone is misleading, as this interpretation is in essence just $\boldsymbol{\mathcal{W}}$ at the first time step. 

The DMD model also consists of a coherent spatial pattern, which we denote $\mathbf{S}$, which is the magnitude of $\boldsymbol{\phi}$,
\begin{equation}\label{eq:DMD-spatial-pattern}
\mathbf{S}_j = |\boldsymbol{\phi}_j| = \sqrt{\boldsymbol{\phi}^{R,2}_j + \boldsymbol{\phi}^{I,2}_j}.  
\end{equation}
$\mathbf{S}$ can be interpreted as the spatial pattern which does not participate in the oscillations of the waveform, i.e., it specifies where the amplitude of $\boldsymbol{\mathcal{W}}$ is the strongest or weakest and is the spatial pattern that does not have a time varying component. Since $\mathbf{S}$ is a strictly positive quantity, any sign changes must be carried by $\boldsymbol{\mathcal{W}}$. The term $2b_j \boldsymbol{S}_j$ describes the full amplitude of the waveform's spatial structure with the factor of two describing the combined conjugate pairs. When fitting spatially stationary data, i.e., with no spatially-varying phase shift, the spatial mode from the phasor notation and the typical DMD interpretation are identical, i.e. $\boldsymbol{\phi}^{R} = |\boldsymbol{\phi}|$ as there is no spatially varying phase shift and thus $\boldsymbol{\phi}^{I}$ must be negligible.

Thus, the DMD model can be re-written as 
\begin{equation}\label{eq:DMD-model-phasor}
\tilde{\mathbf{x}} = \sum_{j=1}^r b_j \mathbf{S}_j \boldsymbol{\mathcal{W}}_j \exp(\mu_jt).
\end{equation}
Taken together, the phasor notation interpretation of the DMD model consists of spatial pattern, $\mathbf{S}$, modulating a spatiotemporal waveform, $\boldsymbol{\mathcal{W}}$, and is scaled by an amplitude, $b$. The exponential growth/decay given by $\mu$ retains its original interpretation. 

The reconstruction of a conjugate pair is
\begin{equation}\label{eq:recon_single-conj-pair}
\tilde{\mathbf{x}}_{j, j+1}(t) = 2b_j \mathbf{S}_j \boldsymbol{\mathcal{W}}_j \exp(\mu_jt)
\end{equation}
where, as in Eq. \ref{eq:added-pairs}, the conjugate pair evenly splits the amplitude of the coherent spatiotemporal mode they describe, as shown by the factor of 2 in Eq. \ref{eq:recon_single-conj-pair}. As long as conjugate pair solutions are found, each pair carries half of the waveform of the reconstructed signal since $\tilde{x}_j$ and  $\tilde{x}_{j+1}$ are in phase with each other (see alsoEq. \ref{eq:added-pairs}).

Thus, the DMD should be interpreted as a collection of \textit{strictly real} \textit{spatiotemporal waveforms} each modulated by a \textit{strictly real and positive} coherent spatial pattern. All of the phasor notation components are real, including the reconstructed data, as long as the DMD solution closely approximates conjugate pairs. If non-conjugate pair solutions are present, the imaginary part of $\tilde{x}$ is not guaranteed to cancel and the physical interpretation loses validity. From this perspective, it may make more sense to label DMD modes according to conjugate pairs instead of rank as this interpretation suggests that we should interpret conjugate pairs only in combination. Similarly, the interpretation of the odd mode (the rank without a conjugate pair) from a DMD model with an odd number of ranks is similarly difficult to reconcile with the phasor notation since this mode will contain an imaginary component even when fitting to strictly real data, which appears physically inconsistent.  

A common strategy to improve a DMD fit is to apply time-delay embeddings \cite{tu2014spectral,brunton2017chaos,arbabiErgodicTheoryDynamic2017}, created by stacking time delayed views of the data as spatial snapshots prior to fitting the DMD model. Fortunately when testing the use of time-delay embeddings, both $\mathbf{S}$ and $\boldsymbol{\mathcal{W}}$ were robustly recovered using $\boldsymbol{\phi}$ from any individual delay. In the PyDMD implementation of the phasor notation we employ the strategy of using the first time delay.


While some of these clarifications and differences between the phasor notation and classical DMD interpretation may seem minor, the phasor notation of DMD can be used to clearly reveal properties of the input data that were difficult to determine previously, as we demonstrate next. These differences become substantial for non-trivial data, especially for windowed DMD approaches applied to more complicated multi-scale as in Sect. \ref{sect:mrCOSTS-phasor-notation}.

\section{Demonstrating the phasor-notation DMD interpretation}\label{sect:uniscale-results}

The utility of the DMD phasor notation is demonstrated using a simple, uniscale toy model that is similar to Tutorial 1 of the PyDMD package \cite{ichinagaPyDMDPythonPackage2024a}. For the uniscale toy model, we create data that are not spatially stationary by summing two functions:
\begin{align}
f_1(x,t) &= \sech(x+3)\cos(2.3t + \frac{x}{10})\\
f_2(x,t) &= 2 \sech(x)\tanh(x)\sin(2.8t + 2.5x).
\end{align}
$f_1$ and $f_2$ are determined on the domain $x=-5...5$ and $t\in[0,4\pi]$. The toy model data, $f_1 + f_2$, is real and composed of 2 spatiotemporal features which oscillate in time with distinct frequencies and spatially-varying phase shifts. A successful DMD model should both faithfully reconstruct the input data and disambiguate the 2 spatiotemporal modes. The spatially-varying phase shift in the sine functions provides two extremes, one with a small shift across the domain ($f_1$) and the other which undergoes multiple oscillations across the domain ($f_2$).  It is important to note that the prescribed phase shift in $f_1$ and $f_2$  are not equivalent to the phase shift from the DMD model Eq. \ref{eq:phase-shift-single-mode} as a consequence of $\mathbf{S}$ being strictly positive (i.e., the sign change from the $\tanh{}$ function will be contained in $\boldsymbol{\mathcal{W}}$)

We define the spatial pattern of each function as the part of the function without the spatially-varying phase shift or the time oscillation
\begin{align}
\hat{f}_1(x) &= |\sech(x+3)|\\
\hat{f}_2(x) &= |2 \sech(x)\tanh(x)|.
\end{align}
A correct definition of the DMD spatial pattern should recover these shapes. For the typical DMD interpretation, the shape should simply be given by the real component of $\boldsymbol{\phi}$, $\boldsymbol{\phi}^{R}$, and for the phasor notation it should be $|\boldsymbol{\phi}|$. To recover the full values of $\hat{f}_1$ and $\hat{f}_2$ from the DMD fit we require the spatial pattern ($\mathbf{S}$), amplitude ($b$), and a factor of 2 from the conjugate pairs. For the typical DMD interpretation this expression is $2b\boldsymbol{\phi}^{R}$ while for the phasor notation this expression is $2b|\boldsymbol{\phi}|$.

\begin{figure}[h]\includegraphics[width=\textwidth]{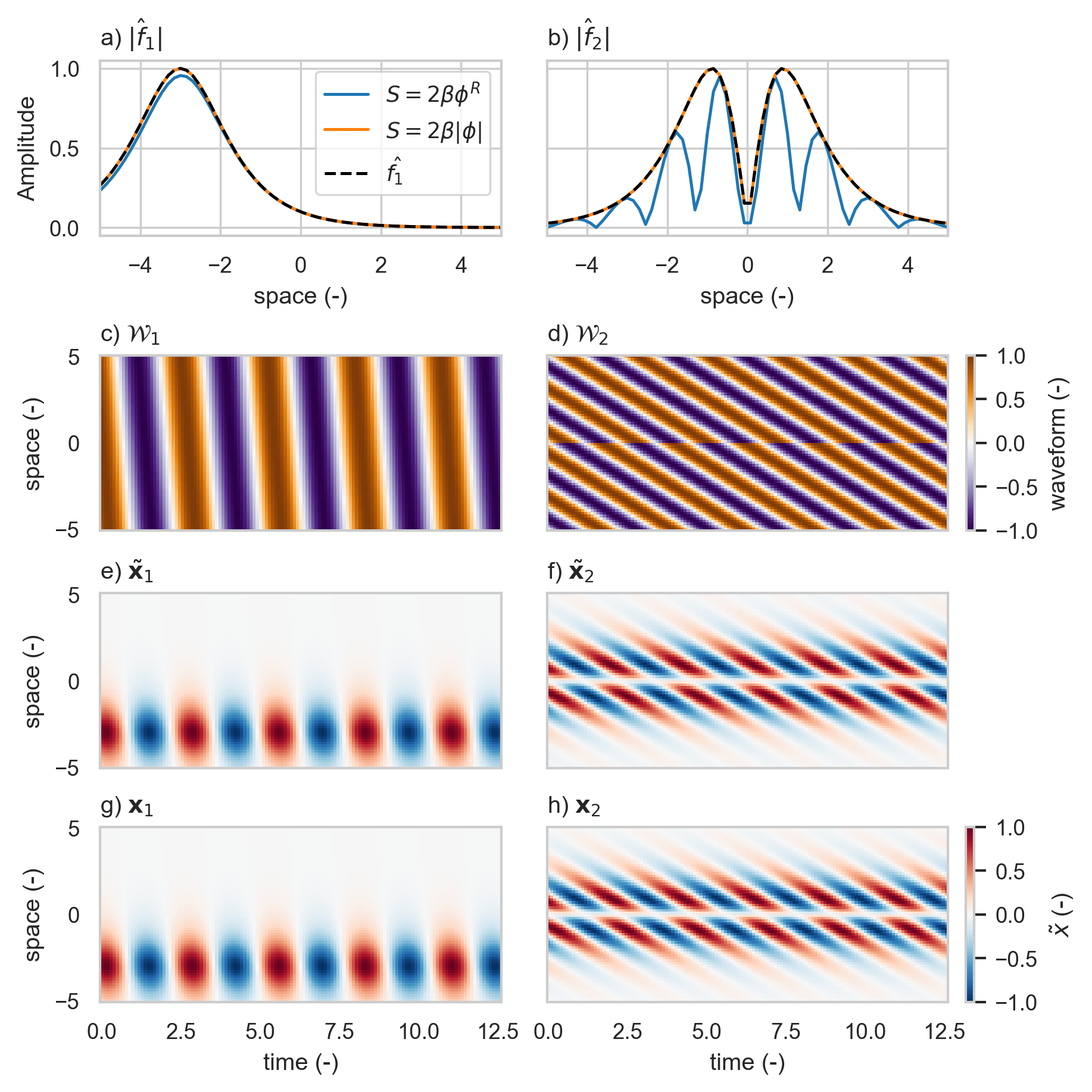}
\caption{A demonstration of $\mathbf{S}$ and $\boldsymbol{\mathcal{W}}$ for the uniscale toy model DMD fit. (a,b) The function amplitudes (dashed black line) with the equivalent spatial patterns from the traditional DMD interpretation (blue line) and the phasor notation interpretation (organge line). (c,d) The spatiotemporal waveforms, $\mathbf{\mathcal{W}}$ (Eq. \ref{eq:DMD-waveform}) are plotted for both modes. (e,f) The reconstruction of each mode in comparison to the (g, h) original input modes are shown. The total relative error of the DMD reconstruction is 1.5 10$^{-7}$.}
\label{fig:uni-scale_spatial-pattern}
\end{figure}

We fit $f_1 + f_2$ with the variable projection solver for DMD \cite{Askham2018} using eigenvalue constraints to force strict conjugate pairs but without a constraint of the real eigenvalues \cite{ichinagaPyDMDPythonPackage2024a}.  The reconstruction of the input signal as well as each individual mode has relative errors on the order of 10$^{-7}$ (Fig. \ref{fig:uni-scale_spatial-pattern}), capturing the input dynamics with high-fidelity.

The spatial patterns, $\hat{f}_1$ and $\hat{f}_2$, are recovered using the phasor notation expression (orange line Fig. \ref{fig:uni-scale_spatial-pattern}a,b) while the traditional DMD interpretation is inconsistent with $\hat{f}_1$ and especially $\hat{f}_2$. For $f_1$, which has a spatial wavelength larger than the size of the spatial domain and thus a smaller spatially-varying phase shift, the typical DMD interpretation largely recovers $\hat{f}_1$ (Fig. \ref{fig:uni-scale_spatial-pattern}a). However, for $f_2$, which has a spatial wavelength shorter than the size of the spatial domain and thus a significant spatially-varying phase shift, the typical DMD interpretation is unable to recover $\hat{f}_2$ (Fig. \ref{fig:uni-scale_spatial-pattern}b).

These results are explainable through the contribution of $\boldsymbol{\phi}^{I}$ to $\tilde{\mathbf{x}}$ through $\mathbf{S}$. For $f_1$ $\boldsymbol{\phi}^{I}$ is small in magnitude and $\boldsymbol{\phi}^{R}$ dominates the results as a consequence of the a low-frequency spatial oscillations. In contrast, for $f_2$ both $\boldsymbol{\phi}^{R}$ and $\boldsymbol{\phi}^{I}$ contribute to $\tilde{\mathbf{x}}$ as a consequence of the high-frequency spatial oscillations. Thus, the spatial patterns are only recoverable using the interpretation provided by the phasor notation, highlighting specifically that $\boldsymbol{\phi}^{I}$ should not be neglected. 

\begin{figure}[h]\includegraphics[width=\textwidth]{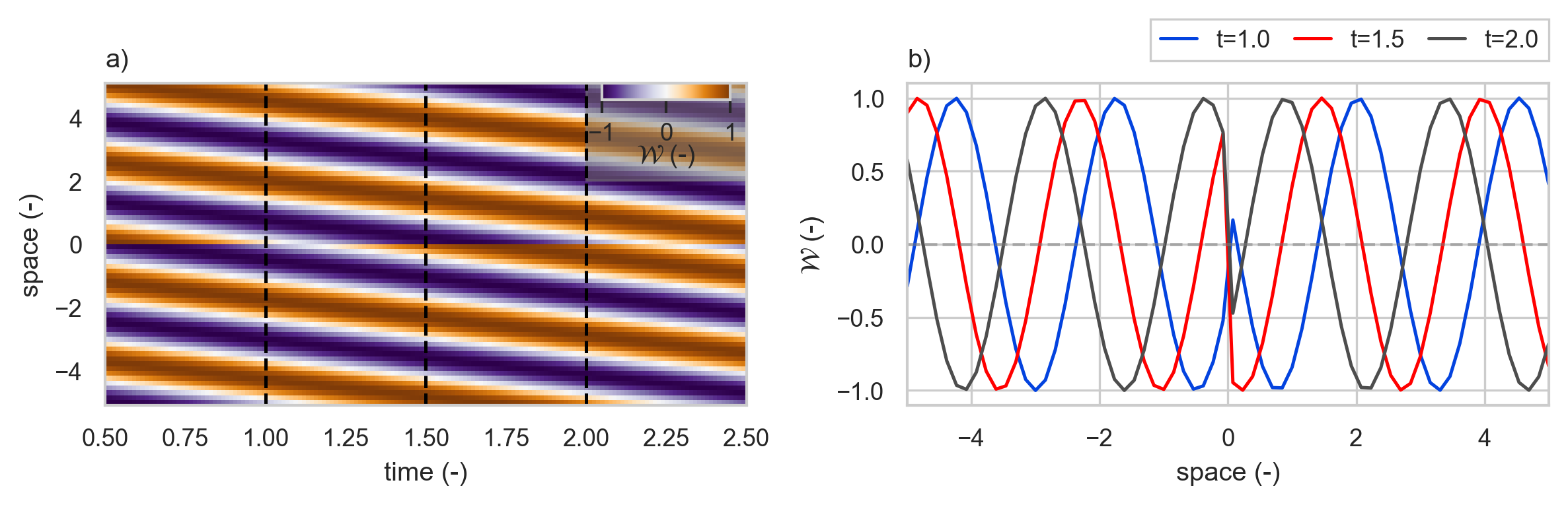}
\caption{a) a zoom in on $\boldsymbol{\mathcal{W}}$ from Fig. \ref{fig:uni-scale_spatial-pattern}d highlight three discrete times ($t$=1, 1.5, and 2, indicated in vertical dashed lines) which are plotted in (b). The color scale in (a) is the same as in Fig. \ref{fig:uni-scale_spatial-pattern}}
\label{fig:uni-scale_waveform-vs-phase-shift}
\end{figure}

In regards to $\boldsymbol{\mathcal{W}}$, the spatially-varying phase shift, $\boldsymbol{\varphi}$, is not a novel component of the typical DMD interpretation. However, one should not interpret $\boldsymbol{\varphi}$ alone, as has been recommended \cite{Taira2020}. As previously stated, such an interpretation is incomplete and can be misleading since $\boldsymbol{\varphi}$ alone is in essence just $\boldsymbol{\mathcal{W}}(t=0)$ and may thus miss important dynamics of the system. 

We demonstrate such a case for $f_2$ in Fig. \ref{fig:uni-scale_waveform-vs-phase-shift}, although qualitatively similar results are found for $f_1$. Here $\boldsymbol{\mathcal{W}}$ is plotted as a function of space at various times. The spatiotemporal waveform of $f_2$ is especially difficult to describe, as the waveform alternates in time between being symmetric and asymmetric in space.  These changes with time make it so that an interpretation at any one time cannot be used to reliably infer the values at another time. The reason $\boldsymbol{\mathcal{W}}$ is not amenable to a simple, intuitive interpretation is due to the angle sum between the components $\boldsymbol{\varphi}$ and $q_jt$ Eq. \ref{eq:DMD-waveform}. One must instead interpret the combined effects of these components as a spatiotemporal wave consisting of a spatially-varying function modulated by an oscillation in time with a single frequency. 

It is commonly asserted that one should examine some combinations of the real and imaginary components of $\boldsymbol{\phi}$ as well as its phase and magnitude to interpret DMD results, e.g., \cite{Taira2020}.  However, since the DMD model is characterized by $\mathbf{S}$ as well as the non-intuitive behavior of the angle sum in $\boldsymbol{\mathcal{W}}$ (Eq. \ref{eq:DMD-phasor-notation}), the real and imaginary components of $\boldsymbol{\phi}$ should not be interpreted -- neither component of $\boldsymbol{\phi}$ describes meaningful physically interpretable results from the DMD model. 

The combination of $\mathbf{S}$ and $\boldsymbol{\mathcal{W}}$ provides an intuitive way of interpreting the original complex functions $f_1$ and $f_2$. For instance, one can mentally multiply $\mathbf{S}$, $b$, and $\boldsymbol{\mathcal{W}}$ to recover each individual function comprising the toy data. It also allows one to examine the coherent spatial pattern directly and accurately. While $\boldsymbol{\mathcal{W}}$ can be more difficult to interpret, it provides a physically-consistent view of the DMD model's waveform. For these reasons, we argue the traditional interpretation of DMD needs to be discarded and replaced by the phasor notation which provides a physically-consistent way to interpret the DMD model.  

\section{Extending the phasor-notation for DMD variants}\label{sect:mrCOSTS-phasor-notation}
A large number of variants for DMD now exist, many of which alter the relationship between the DMD model and the fitted data. For these cases, one must take care to adjust the phasor notation expressions to be consistent with the DMD variant. For instance, some DMD variants fit the DMD model to time subsamples, resulting in a collection of DMD models capable of diagnosing complex behaviors such as transience and invariances. Examples include the non-stationary DMD (nsDMD) \cite{ferreNonStationaryDynamicMode2023c}, multi-resolution DMD \cite{Kutz2016_dmd, Dylewsky2019_dmd}, and mrCOSTS \cite{lapoMethodUnsupervisedLearning2025}. In these cases one must represent the summed effects of the DMD components from overlapping time windows. Here, we develop the phasor notation for mrCOSTS, but each method has their own idiosyncrasies requiring specific development of the phasor notation. This example is intended to demonstrate a strategy for deriving similar expressions for other DMD variants. 

At a given decomposition level, mrCOSTS fits a DMD model to data using sliding, overlapping windows in time. It then separates out an unresolved low-frequency band from the resolved high-frequency bands. The low frequency band is then passed to the next level with a larger window size which is capable of resolving more of the low-frequency component of the signal. This step is repeated to a desired decomposition level at which point a global scale separation is performed to capture dynamics leaked between levels, yielding a hierarchical decomposition of bands describing a relatively narrow range of frequencies \cite{lapoMethodUnsupervisedLearning2025}.

Bands are denoted using the index $p$, which is comprised by a subset of the ranks ($j$), time windows ($k$), and decomposition levels ($\ell$) of DMD modes that exist within the frequency band, defined as $(j, k, \ell) \in p$. The reconstruction of a band in this situation can be notated as
\begin{equation}\label{eq:mrCOSTS-summed}
\tilde{\mathbf{x}}_p(t)=\sum_{(j, k, \ell) \in p} b_{j, k, \ell} e^{\lambda_{j, k, \ell}t} \boldsymbol{\phi}_{j, k, \ell}
\end{equation}
which can be re-written in phasor notation as
\begin{equation}\label{eq:DMD-summed-phasor}
\tilde{\mathbf{x}}_p(t)=\sum_{(j, k, \ell) \in p} b_{j, k, \ell} e^{\mu_{j, k, \ell}t} \textbf{S}_{j, k, \ell} \boldsymbol{\mathcal{W}}_{j, k, \ell}.
\end{equation}
The complexity of Eq. \ref{eq:DMD-summed-phasor} is that it includes the superposition of spatiotemporal waveforms and spatial patterns including modulations of $\boldsymbol{\mathcal{W}}$ between time windows, which can result in beat frequencies, as well as variations in $\mu_{j, k, \ell}$. Phasor notation expressions for mrCOSTS need to be able to represent these complexities.

\subsection{mrCOSTS phasor notation summed terms} 
A reasonable expression for the band amplitude, $\beta_p(t)$, is
\begin{equation}\label{eq:summed-amplitude}
\beta_p(t) = \sum_{(j, k, \ell)\in p} e^{\mu_{j, k, \ell}t} b_{j, k, \ell}
\end{equation}
which includes the non-oscillatory time component, $e^{\mu t}$, in order to take into account windows with strong non-linearities. To represent the summed effects for $\mathbf{S}$ and $\boldsymbol{\mathcal{W}}$ we use a weighted average, where the band amplitude is used to weight the summed expressions. The spatial pattern of the summed bands, $\mathbf{S}_p$, is the weighted average of all the terms modulating the spatial pattern's amplitude
\begin{equation}\label{eq:summed-spatial-pattern}
\mathbf{S}_p(x, t) = \frac{\sum_{(j, k, \ell)\in p} e^{\mu_{j, k, \ell}t} b_{j, k, \ell} \sqrt{\boldsymbol{\phi}^{R, 2}_{j, k, \ell} + \boldsymbol{\phi}^{I,2}_{j, k, \ell}}}{\beta_p(t)}.
\end{equation}
The weighted average is necessary so that modes with a small value for $b$ or negative $\mu$ are recognized as contributing less to the resulting $\mathbf{S}_p$ than modes with a large $b$ or a positive $\mu$.

We follow a similar logic to determine the summed $\boldsymbol{\mathcal{W}}_p$. We cannot simply sum the waveforms themselves, as a waveform with a small amplitude will contribute less to the waveform of band $p$ than one with a large amplitude. Thus, we choose to also represent $\boldsymbol{\mathcal{W}}$ as a weighted average,
\begin{equation}\label{eq:summed-waveform}
\boldsymbol{\mathcal{W}}_p(x, t) = \frac{\sum_{(j, k, \ell)\in p} b_{j, k, \ell} e^{\mu_{j, k, \ell}t} \cos(\omega{j, k, \ell}t + \boldsymbol{\varphi}_{j, k, \ell})}{\beta_p(t)}.
\end{equation}

The weighted average expressions allow for capturing complex beat frequencies that can emerge, even when summing over fairly broad frequency bands. Consider two overlapping time windows, $k=1$ and $k=2$ where $k=1$ proceeds $k=2$. Window $k=1$ has a DMD mode with large $b$ but also has a large and negative $\mu$; the mode rapidly decays across the window. Window $k=2$, which overlaps with $k=1$, has a DMD mode with a smaller $b$ but with $\mu$ near to zero. In the time range where the windows overlap, the rapidly decaying mode in window $k=1$ contributes less to $\boldsymbol{\mathcal{W}}_p(x, t)$ and $\mathbf{S}_p(x,t)$ as result of the temporal decay than the mode without rapid decay in window $k=2$.

These expressions are only approximations of the phasor notation terms. Specifically, through the weighted average, we assume that the individual components vary slowly between windows. Thus, it is important to note that $\tilde{\mathbf{x}}$ can only be approximated when using these summed terms. Further, if a decomposition includes a rapid and large variation in one of the terms in $\beta_p$, then the weighted average approach loses validity. 

\subsection{Multi-scale toy model}
To demonstrate the mrCOSTS phasor notation terms, a multi-scale toy model is implemented that is similar to the tutorial for the PyDMD implementation of mrCOSTS \cite{ichinagaPyDMDPythonPackage2024a} as well as the toy model presented in \cite{Dylewsky2019_dmd} for mrDMD. The multi-scale toy model consists of two processes with distinct, but overlapping time scales: the Fitzhugh-Naguma model, which has a spiking behavior, and the Unforced Duffing Oscillator, which is a non-linear oscillator. A transient, translating wave packet with a high frequency is then added to the oscillators.

\begin{figure}[h]\includegraphics[width=\textwidth]{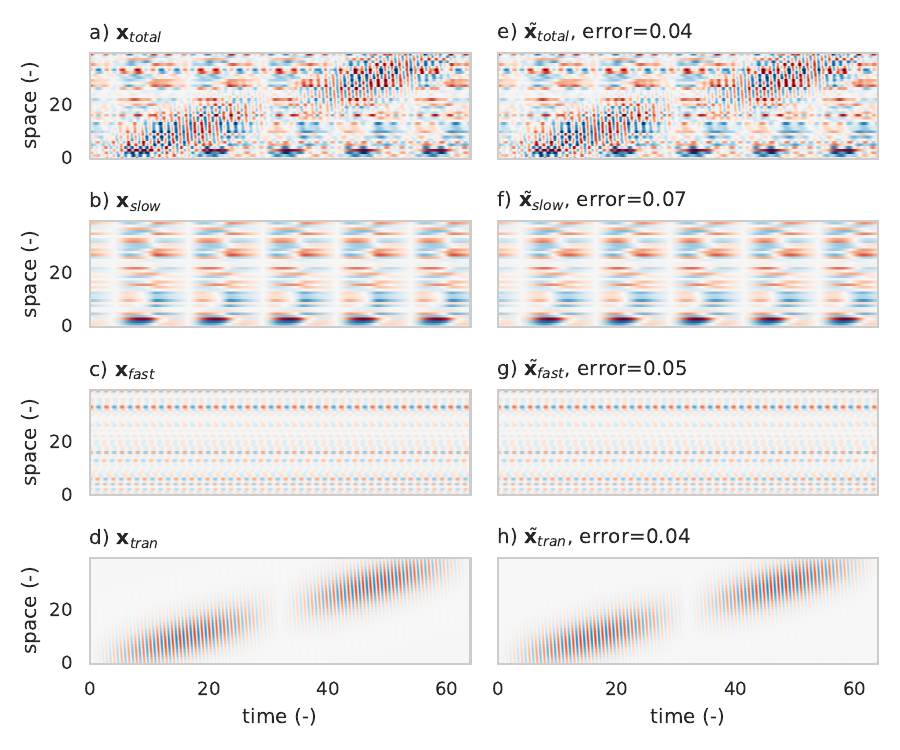}
\caption{(a) The full multi-scale toy model and (b-d) the individual components are shown. (e) The corresponding mrCOSTS reconstruction of the full signal as well as (f-h) the individual components of the mulit-scale toy model are plotted. The error of the reconstructions is specified in the subtitles.}
\label{fig:multi-scale_recon}
\end{figure}

The Fitzhugh-Naguma model is the slow mode (Fig. \ref{fig:multi-scale_recon}b) described by
\begin{align}\label{eq:slow-mode}
    \dot{v}& = v - \frac{1}{3}v^3 - w + 0.65 \\
    \dot{w} &= \frac{1}{\tau_1}(v + 0.7 - 0.8w).
\end{align}
The Unforced Duffing Oscillator is the fast mode (Fig. \ref{fig:multi-scale_recon}c) described by 
\begin{align}\label{eq:fast-mode}
    \dot{p}& = q \\
    \dot{q}& = -\frac{1}{\tau_2}(p + p^3).
\end{align}
We impose two distinct time scales, $\tau_1= 2$ and $\tau_2=0.2$, but due to the non-linear nature of the models the slower time scale of Eq. \ref{eq:slow-mode} can approach the time scale of the faster time scale in Eq. \ref{eq:fast-mode}. A combined signal is generated by stacking $[v, w]$ and $[ p,q]$ ten times, representing an arbitrary spatial dimension $\mathbf{x}$ with a length of 40. These two models are then linearly mixed using a random orthogonal matrix $O$, yielding

\begin{equation}
\mathbf{x}_{slow}+\mathbf{x}_{fast}=
    \begin{bmatrix}
        u&v& \cdots & u &v&p&q&\cdots&p&q
    \end{bmatrix}\times O.
\end{equation}

Finally, a transient, translating high-frequency wave packet with a sinusoidal waveform (Fig. \ref{fig:multi-scale_recon}d) is created following
\begin{align}\label{eq:transient}
\mathbf{x}_{tran}&= 2 \sin(10 t + \frac{x}{2 \pi}) \sigma_x(x,t)\sigma_t(t)\\
\boldsymbol{\sigma}_x(x,t)&= \exp\left(-\frac{(x - \frac{x_0t}{T})^2}{(0.2x_0)^2}\right)\\
\sigma_t(t) &= |\sin\left(\frac{2\pi t}{T}\right)|
\end{align}
where $x_0$ is the maximum value of the x domain, and $T$ is the length of the time domain. $\sigma_x(x,t)$ is a Gaussian distribution that travels along the space dimension with time and describes the spatial pattern of the wave packet. We prescribe the transient behavior of the wave packet using $\sigma_t(t)$ yielding one maxima in the first half of the time domain and another in the second half of the time domain. These data are generated on the domain $\mathbf{x} = 0....40$ and $t=0....64$. The sum of all three signals yields the complete multi-scale model, $\mathbf{x}_{total}$ (Fig. \ref{fig:multi-scale_recon}a). The translating and transient behavior of $\mathbf{x}_{tran}$ and the varying time scales $\mathbf{x}_{slow}$ and $\mathbf{x}_{fast}$ cannot be diagnosed using a typical DMD model.

\subsection{Multi-scale phasor notation interpretation} The multi-scale toy model was decomposed using mrCOSTS with three decomposition levels (window lengths of [60, 120, 480] with a window slide of 10\% between windows). Each level was fit using a rank of six with an eigenvalue constraint requiring conjugate pairs. 

\begin{figure}[h]\includegraphics[width=\textwidth]{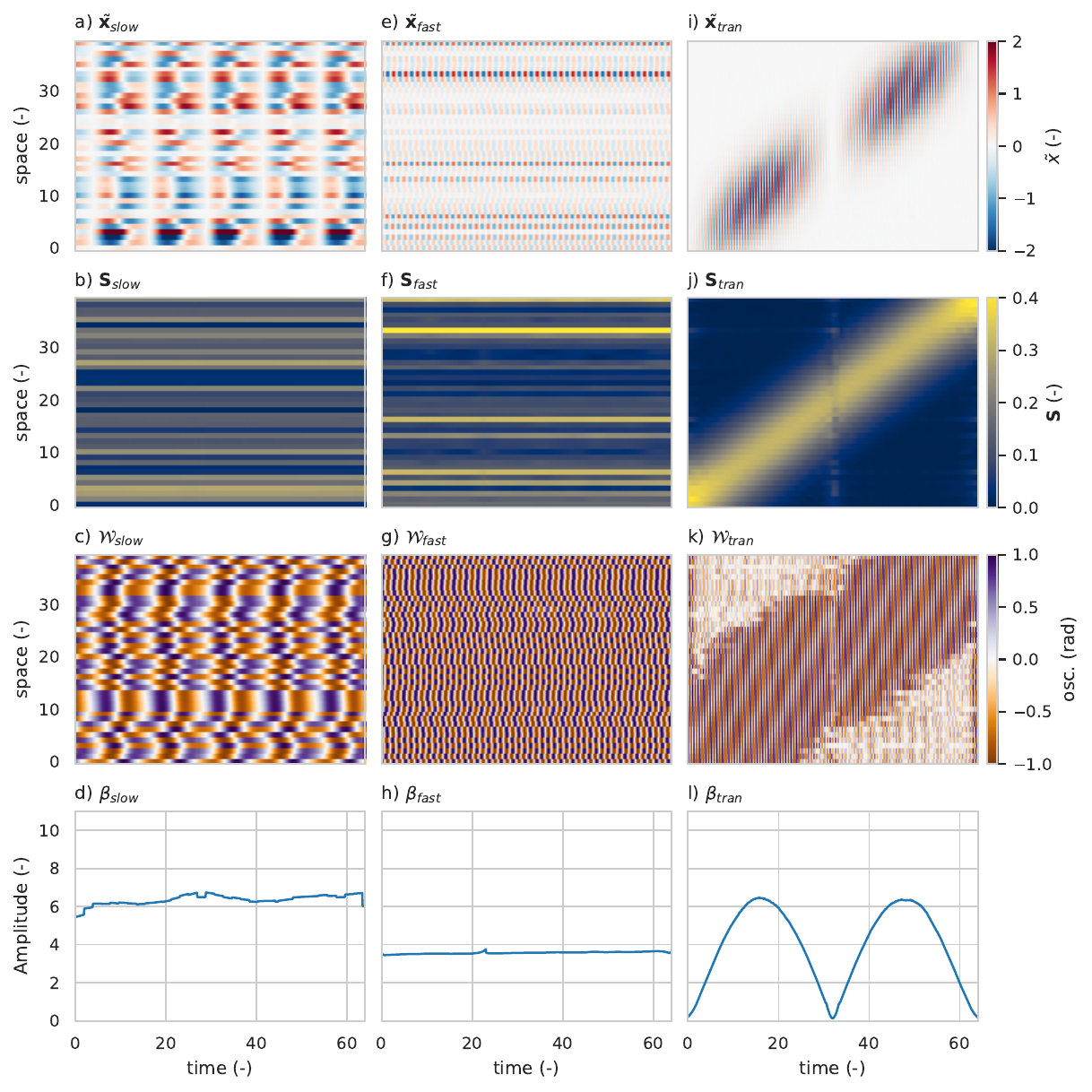}
\caption{Summed phasor notation components and reconstruction for the mrCOSTS decomposition of the multi-scale toy model. (a-d) The left column shows the slow modes, (e-h) the middle column the fast modes, and (i-l) the right column the transient mode. The first row (a, e, i) shows the reconstruction of each mode, the second row  (b, f, j) shows $\mathbf{\mathcal{S}}_p$ for each mode, the third row  (c, g, k) shows $\boldsymbol{\mathcal{W}}_p$ for each mode, and the bottom row  (d, h, l) shows $\beta$ for each mode.}
\label{fig:multi-scale_components}
\end{figure}

The mrCOSTS reconstruction has a high-fidelity given the multi-scale nature of the data with a relative error of 4.5\% (Fig. \ref{fig:multi-scale_recon}e). Each component of the multi-scale toy model system is successfully identified and decomposed (Fig. \ref{fig:multi-scale_recon}b-d vs Fig. \ref{fig:multi-scale_recon}f-h). Generally, mrCOSTS reconstructions have larger errors than those found for DMD fits to uniscale data, but this is expected and the constituent components are recovered with sufficient fidelity, with errors between 3.6\% (transient mode) to 6.8\% (slow mode; Fig. \ref{fig:multi-scale_recon}). 

The spatiotemporal variations are attributable to the summed waveform, $\boldsymbol{\mathcal{W}}_p$ (Fig. \ref{fig:multi-scale_components}c,g,k), which describes the oscillations with no regard to spatial pattern or time varying amplitude. Consequently, each mode's $\boldsymbol{\mathcal{W}}_p$ scales between -1 and 1 and the amplitude of the oscillatory pattern in each $\boldsymbol{\mathcal{W}}_p$ is unaffected by the mode's amplitude or the specific spatial structure of the function, demonstrating  that the weighted averaging in Eq. \ref{eq:summed-waveform} behaves as expected. 

For $\boldsymbol{\mathcal{W}}_{tran}$ (Fig. \ref{fig:multi-scale_components}k), artifacts are present in the regions where the amplitude of $\mathbf{x}_{tran}$ is sufficiently small that mrCOSTS can no longer resolve it (top left and bottom right regions of Fig. \ref{fig:multi-scale_components}k). In the regions where $\mathbf{x}_{tran}$ is of sufficient amplitude, $\boldsymbol{\mathcal{W}}_{tran}$ reveals a pure sine wave with a spatially varying phase shift. The imprint of the spatial translation is also apparent in $\boldsymbol{\mathcal{W}}_{tran}$, but the spatial structure, namely the Gaussian distribution, is not. Similarly, the time varying amplitude of the wave packet is not present in $\boldsymbol{\mathcal{W}}_{tran}$.  These results for $\boldsymbol{\mathcal{W}}_{tran}$ are expected as the Gaussian distribution of the wave packet should be described instead by $\mathbf{S}_{tran}$ and the time varying amplitude should be described by $\beta_{tran}$. Instead, $\boldsymbol{\mathcal{W}}_{tran}$ purely describes the mode's oscillations agnostic to the spatial pattern or time varying amplitude.

A similar interpretation is made for $\boldsymbol{\mathcal{W}}_{slow}$ and $\boldsymbol{\mathcal{W}}_{fast}$ as for $\boldsymbol{\mathcal{W}}_{tran}$. However, these modes are present with sufficient amplitude throughout the entire domain and thus do not have the artifacts present in $\boldsymbol{\mathcal{W}}_{tran}$. In $\boldsymbol{\mathcal{W}}_{fast}$ and $\boldsymbol{\mathcal{W}}_{slow}$, the oscillatory patterns are not pure sine waves, but capture oscillations with time varying frequencies in a relatively narrow band ($\boldsymbol{\mathcal{W}}_{fast}$) and a fairly broad band of frequencies ($\boldsymbol{\mathcal{W}}_{slow}$). In both cases, the oscillatory pattern of the diagnosed modes are clearly captured (e.g., compare Fig. \ref{fig:multi-scale_components}a to Fig. \ref{fig:multi-scale_components}c).

$\mathbf{S}_p$ describes how the spatial pattern of a band varies in time with no regard to variations of the band's amplitude with time, which are instead described by $\beta_p$.  For the slow and fast modes, both $\mathbf{S}_p$ (Fig. \ref{fig:multi-scale_components}b,f) and $\beta_p$ (Fig. \ref{fig:multi-scale_components}c,g), are effectively static in time. This result is expected as the physics of the slow and fast modes were prescribed to be constant in time. $\beta_{slow}$ does contain some artifacts that are a consequence of the mrCOSTS fits degrading at the edges of the time domain and for the largest decomposition level.

In contrast to the fast and slow modes, $\mathbf{x}_{tran}$ was prescribed to translate in space with a time varying amplitude. The translation is described by $\mathbf{S}_{tran}$ (\ref{fig:multi-scale_components}j), which shows the translation of the spatial pattern independent of the oscillations of the mode's amplitude. The mode's amplitude variations are well-described by  $\beta_{tran}$ (Fig. \ref{fig:multi-scale_components}l), which captures the sinusoidal amplitude pattern we imposed in Eq. \ref{eq:transient}. There is an artifact in the middle of the time domain where $\mathbf{x}_{tran}$ becomes small enough that mrCOSTS can no longer describe it. Otherwise, the prescribed behaviors in Eq. \ref{eq:transient} are well-described by the mrCOSTS phasor notation terms, demonstrating how these terms can be used to recover a physically-consistent interpretation of the diagnostic. 

\section{Conclusions and Outlook}
\label{sec:conclusions}

A definition of DMD as a model, rather than an algorithm, was presented. In the DMD model, modes consist of strict conjugate pairs, each of which carries half of the mode's waveform, with the imaginary components of each pair perfectly canceling each other. When a DMD model is fit to real data without using strict conjugate pairs, the solution will not perfectly cancel out the imaginary components creating an inconsistent physical interpretation. 

Using the existence of strict conjugate pairs, the phasor notation of the DMD model was presented. The phasor notation allows a physically consistent interpretation of the DMD model that was previously missing. In this new view of the DMD model, a DMD mode consists of a spatial pattern that is the magnitude of the complex variable, $\boldsymbol{\phi}$ and a waveform described by the imaginary component of the mode's eigenvalue, which describes the oscillatory frequency of the waveform, and the ratio of the imaginary and real parts of $\boldsymbol{\phi}$, which describes the spatially varying phase shift. The phase shift at a single point in time or the real or imaginary parts of $\boldsymbol{\phi}$ are generally not physically interpretable. 

This new interpretation of the DMD model enables a simple physical interpretation through a subtle, but important, re-interpretation of the terms contained in it. This new interpretation was demonstrated on a simple uni-scale toy model for the classic DMD model. For variants of DMD that alter the relationship between the model and the underlying data, different phasor notation interpretations are necessary. We demonstrated one such interpretation for mrCOSTS, a hierarchical, multi-resolution variant of DMD, presenting expressions that closely approximate the phasor notation terms for this DMD variant. 

Consequently, the phasor notation of the DMD model addresses the interpretability gap. The interpretability gap is a term we introduce which refers to the difficulty in bridging between the mathematical properties of a method and the physical interpretation the method can provide.  Closing interpretability gaps has been a powerful way of integrating a new method into the physical sciences, such as the case of presenting continuous wavelet transforms (CWT) to the geophysical community by \cite{Torrence1998}. We argue this gap is of more general importance and that care should be taken to verify and include a physical understanding of a method, as this is a critical step in the method's implementation. 

Bridging the interpretability gap can benefit from a closer integration of the expertise in both the mathematical and physical sciences, as this relationship can often be more unidirectional rather than a dialogue. As concrete example we draw attention to the difference between how DMD was described previously and the physical interpretation provided in this work or the difference between how CWTs were used in a more qualitative sense prior to the work of Torrence and Compo \cite{Torrence1998}. As such, we argue the interpretability gap should be more explicitly considered in future work and is not a niche concern.

\section*{Acknowledgments}
This research was funded in part by the Austrian Science Fund (FWF) [10.55776/ESP214]. This work was supported in part by the US National Science Foundation (NSF) AI Institute for Dynamical Systems (dynamicsai.org), grant 2112085. JNK further acknowledges support from the Air Force Office of Scientific Research  (FA9550-24-1-0141).

\bibliographystyle{unsrtsiam}
\bibliography{library}
\end{document}